\documentclass[11pt]{article}
\usepackage{mathdots}

\usepackage{fullpage}
\usepackage{authblk}
\usepackage{natbib}
\usepackage{algorithm}
\usepackage{algpseudocode}
\usepackage{amsmath,amsthm,amsfonts}
\usepackage{amsmath}
\usepackage[subnum]{cases}
\usepackage{tcolorbox}
\usepackage{mathrsfs}
\usepackage{amssymb}
\usepackage{color}
\usepackage{mathrsfs}
\usepackage{empheq}
\usepackage{enumitem}
\usepackage{bm}
\usepackage{multirow}
\usepackage{booktabs}
\usepackage{makecell}
\usepackage{graphicx}
\usepackage{subcaption}
\usepackage{comment}
\usepackage{cases}
\usepackage{appendix}
\usepackage{tikz}
\usetikzlibrary{arrows,shapes}
\usepackage[colorlinks= true, linkcolor=brown, citecolor=blue, urlcolor=black]{hyperref}
\usepackage{cleveref}
\usepackage{marginnote}
\usepackage{diagbox} 

\numberwithin{equation}{section}

\theoremstyle{definition}

\newtheorem{theorem}{Theorem}[section]

\newtheorem{remark}[theorem]{Remark}

\newtheorem{defn}[theorem]{Definition}

\DeclareMathOperator*{\argmin}{arg\,min}

\usepackage{booktabs}
\usepackage{tikz}
\usepackage{pgfplots}
\usetikzlibrary{calc,intersections}

\usetikzlibrary{shapes.geometric, arrows, positioning}

\tikzset{
	mybox/.style  = {draw, rectangle, minimum width=4cm, minimum height=0.8cm, text centered, text width=4.4cm,   
		font=\normalsize},
	box/.style  = {draw, rectangle, minimum width=2.0cm, minimum height=0.6cm, text centered, text width=3.0cm,   
		font=\normalsize},
	myarrow/.style = {line width=0.2pt, draw=black, -triangle 60, postaction={draw, line width=0.2pt, shorten >=10pt,-}}
}

\tikzstyle{arrow} = [->, >=stealth, -triangle 60]

\allowdisplaybreaks

\makeatletter
\newcommand{\leqnomode}{\tagsleft@true}
\newcommand{\reqnomode}{\tagsleft@false}
\makeatother

\begin{document}

\title{Accelerated Implicit GDA Schemes: Theoretical Guarantees and Application to Proximal Augmented Lagrangian Methods}
\author[2,3]{Jiaqi Liu}
\author[1,2]{Bin Shi\thanks{Corresponding author: \url{binshi@fudan.edu.cn} } }
\affil[1]{Center for Mathematics and Interdisciplinary Sciences, Fudan University, Shanghai 200433, China}
\affil[2]{Shanghai Institute for Mathematics and Interdisciplinary Sciences (SIMIS), Shanghai 200433, China}
\affil[3]{Research Institute of Intelligent Complex Systems, Fudan University, Shanghai 200433, China}
\date\today

\maketitle

\begin{abstract}
Convex optimization problems with linear equality constraints arise ubiquitously in scientific computing, machine learning, and control theory. While classical Krylov subspace methods are highly effective, they often rely on problem-specific preconditioners and incur significant memory overhead. Conversely, gradient-based methods, such as the augmented Lagrangian method (\texttt{ALM}), avoid these limitations but frequently suffer from slow convergence, particularly in their outer iterations. Developing accelerated outer-iteration schemes, therefore, remains a critical research objective. In this study, we demonstrate that incorporating a proximal operation into the augmented Lagrangian framework yields the proximal \texttt{ALM}, where the outer iteration is equivalent to an implicit gradient descent-ascent (\texttt{GDA}) scheme. We further establish that this equivalence extends naturally to the setting of variable step sizes. Through Lyapunov analysis, we show that the underlying potential function must be shifted from the conventional objective gap to a variational inequality measure, signaling a shift in perspective from pure convex optimization to minimax optimization.   Motivated by these observations, we first develop an implicit GDA scheme with variable step sizes based on a continuous-time ODE framework, which achieves an $o(1/k)$ last-iterate convergence rate for both the primal-dual objective gap and the gradient norm. Building upon a second-order ODE framework, we then propose a family of Nesterov-type implicit \texttt{GDA} schemes parameterized by $r \geq 0$, which achieves an $o(1/k^{r+1})$ last-iterate convergence rate for the primal-dual objective gap. Furthermore, specializing the second-order ODE formulation to the case $r=0$, we derive a corresponding explicit \texttt{GDA} scheme and prove an $o(1/k)$ last-iterate convergence rate for the primal-dual objective gap. Finally, we present several numerical experiments to validate these theoretical results and demonstrate the effectiveness of the proposed methods.
\end{abstract}

%

\section{Introduction}
\label{sec: intro}

Convex optimization with equality constraints forms a cornerstone across scientific computing, machine learning, and control theory. A primary application arises in saddle-point systems governing fluid dynamics and electromagnetism, such as the steady-state Stokes and Maxwell equations, where equality constraints enforce physical conservation laws or divergence-free conditions~\citep{benzi2005numerical}. This framework is also central to inverse problems and variational data assimilation, where a regularized misfit function is minimized subject to the underlying partial differential equation (PDE)-constrained dynamics~\citep{le1986variational}. Furthermore, in systems and control, the classical linear quadratic regulator (LQR) problem exemplifies this framework, optimizing state and control trajectories under linear state-space transitions~\citep{anderson2007optimal}.

To motivate our study, we consider the canonical convex optimization problem with linear equality constraints:
\begin{equation}
\label{eqn: cvx-equal-constraint}
\left\{ \begin{aligned}
         & \min f(x) = \frac12 \| Ax - b \|^2, \\
         & \mathrm{s.t.} \quad  Ex = h.
        \end{aligned} \right.
\end{equation}
The problem~\eqref{eqn: cvx-equal-constraint} can be reformulated via the associated Lagrangian function:
\begin{equation}
\label{eqn: lagrangian-multiplier}
L(x; y) = f(x) + \langle y, Ex - h \rangle, 
\end{equation}
where $y$ denotes the vector of Lagrange multipliers. The solution to the constrained optimization problem~\eqref{eqn: cvx-equal-constraint} corresponds to the saddle point of the Lagrangian function~\eqref{eqn: lagrangian-multiplier}, which is equivalent to solving the following block linear system, known as the Karush-Kuhn-Tucker (KKT) conditions:
\begin{equation}
\label{eqn: matrix-form-saddle}
\begin{pmatrix}
A^{T}A &  E^{\top} \\
E         &   0 
\end{pmatrix} \begin{pmatrix} x \\ y \end{pmatrix} = \begin{pmatrix} A^{\top}b \\ g \end{pmatrix}.
\end{equation}
Conventionally, this large-scale saddle-point system~\eqref{eqn: matrix-form-saddle} is solved via Krylov subspace methods~\citep{benzi2005numerical}. Among these approaches, the preconditioned generalized minimal residual method (GMRES)~\citep{benzi2004preconditioner}, building on the foundational work of~\citet{saad1986gmres}, is one of the most widely used techniques. When the system can be formulated within a symmetric framework, the minimal residual method (MINRES)~\citep{paige1975solution} is often preferred. Nevertheless, the preconditioned GMRES suffers from increasing memory and computational overhead as iterations progress, and its practical performance depends heavily on the availability of high-quality preconditioners tailored to the structure of the underlying problem. As a result, there is no universally effective solver for the large-scale saddle problem~\eqref{eqn: matrix-form-saddle}; instead, efficient algorithms must typically be developed on a case-by-case basis. Furthermore, when the coefficient matrix is nonnormal, spectral analysis alone is insufficient to accurately predict convergence behavior. In such settings, the analysis of pseudospectra becomes essential, constituting an active area of contemporary research~\citep{trefethen2020spectra}.


From a variational perspective, the augmented Lagrangian is defined as:
\begin{equation}
\label{eqn: aug-lagrangian}
L_s = f(x) + \langle y, Ex - h \rangle + \frac{s}{2} \| Ex - h \|^2,
\end{equation}
which exactly possesses the same saddle point as the Lagrangian function~\eqref{eqn: lagrangian-multiplier}. Exploiting this property, the augmented Lagrangian method (\texttt{ALM}), originally proposed by~\citet{glowinski1975sur}, updates the variables via the following iterations:
\begin{subequations}
\label{eqn: aug-lag-method}
\begin{empheq}[left=\empheqlbrace]{align}
         & x_{k} = \argmin_{x \in \mathbb{R}^n} \left\{ f(x) + \frac{s}{2} \left \| Ex - h + \frac{y_{k-1}}{s} \right \|^2 \right\},      \label{eqn: alm-x}   \\
         & y_{k} = y_{k-1} + s \left( E x_{k} - h \right),                                                                                                       \label{eqn: alm-y} 
\end{empheq}    
\end{subequations}
which serves as an efficient numerical optimization method for solving the constrained convex problem~\eqref{eqn: cvx-equal-constraint}. This framework was later popularized by~\citet{boyd2010distributed} for distributed optimization, including applications such as the least absolute shrinkage and selection operator (Lasso). Structurally, the~\texttt{ALM} iteration~\eqref{eqn: aug-lag-method}  operates on two coupled levels:  the inner iteration and the outer iteration. The inner iteration is embedded within the implicit $\mathrm{argmin}$ operation in~\eqref{eqn: alm-x}, which involves solving a symmetric positive semidefinite quadratic optimization problem.  Such subproblems can often be solved efficiently using classical iterative solvers, such as the conjugate gradient method or the successive over-relaxation method. The outer iteration refers to the joint update scheme from $k$ to $k+1$, encompassing both the primal minimization~\eqref{eqn: alm-x} and the dual gradient ascent~\eqref{eqn: alm-y}. The convergence rates of this outer sequence were first established by~\citet{he20121, he2015non}. Most recently, ~\citet{li2024understanding} utilized Lyapunov analysis alongside a high-resolution ordinary differential equation (ODE) framework to investigate the underlying mechanisms of this outer convergence behavior, employing the discrete Lyapunov function:
\begin{equation}
\label{eqn: lyapunov-alm}
\mathcal{E}(k) = \frac{1}{2} \| y_k - y^{\star} \|^2.
\end{equation}
Specifically, for the ergodic (running) averages, defined as
\[
\overline{x}_{k} = \frac{1}{k+1} \sum_{i=0}^{k}x_{i+1}, \quad \mathrm{and} \quad \overline{y}_{k} = \frac{1}{k+1} \sum_{i=0}^{k}y_{i+1},
\]
the convergence rate is bounded by
\begin{equation}
\label{eqn: rate-alm-ave}
L(\overline{x}_k; y^{\star}) - L(x^{\star}; \overline{y}_k) \leq O\left( \frac{1}{k} \right).
\end{equation}
Further analysis of the numerical error indicates that the constraint violation for the last iterate only achieves a rate of:
\begin{equation}
\label{eqn: rate-alm-last}
\| E x_{k} - h  \|^2  \leq O\left( \frac{1}{k} \right).
\end{equation}

Both convergence rates,~\eqref{eqn: rate-alm-ave} and~\eqref{eqn: rate-alm-last}, are unsatisfactorily slow, which motivates the development of accelerated schemes. Inspired by the emerging Lyapunov analysis and the high-resolution ODE framework for gradient-based optimization methods~\citep{shi2022understanding}, we propose several accelerated schemes in this study based on the proximal augmented Lagrangian method 
\begin{subequations}
\label{eqn: prox-aug-lag-method}
\begin{empheq}[left=\empheqlbrace]{align}
         & x_{k} = \argmin_{x \in \mathbb{R}^n} \left\{ f(x) + \frac{s}{2} \left\| Ex - h + \frac{y_{k-1}}{s} \right \|^2 + \frac{1}{2s} \left\| x - x_{k-1} \right \|^2 \right\},          \label{eqn: prox-alm-x}  \\
         & y_{k} = y_{k-1} + s \left( E x_{k} - h \right).                                                                                                                                                                                \label{eqn: prox-alm-y}
\end{empheq}    
\end{subequations}

\subsection{An implicit GDA viewpoint with variable step-size extension}
\label{subsec: gda}

We begin by characterizing the implicit minimization step~\eqref{eqn: prox-alm-x}, which satisfies the first-order optimality condition:
\[
\nabla f(x_{k}) + sE^{\top} \left( Ex_{k} - h + \frac{y_{k-1}}{s} \right) + \frac{x_{k} - x_{k-1}}{s} = 0.
\]
Substituting the dual update~\eqref{eqn: prox-alm-y} into this relation, we obtain:
\begin{equation}
\label{eqn: prox-alm-x-expan-2}
\nabla f(x_{k}) + E^{\top} y_{k}+  \frac{x_{k} - x_{k-1}}{s} = 0. 
\end{equation}
Recalling the Lagrangian defined in~\eqref{eqn: lagrangian-multiplier}, we can reformulate equations~\eqref{eqn: prox-alm-x-expan-2} and~\eqref{eqn: prox-alm-y} as the following implicit gradient descent-ascent (\texttt{GDA}) iteration: 
\begin{subequations}
\label{eqn: implicit-gda}
\begin{empheq}[left=\empheqlbrace]{align}
         & x_{k} = x_{k-1} -  s  \nabla_x L(x_{k},  y_{k}),                      \label{eqn: implicit-gda-x}   \\
         & y_{k} = y_{k-1} + s  \nabla_y L(x_{k},  y_{k}).                       \label{eqn: implicit-gda-y}
\end{empheq}    
\end{subequations}
In other words,  the outer iteration of the proximal~\texttt{ALM}~\eqref{eqn: prox-aug-lag-method} can be viewed equivalently as the implicit~\texttt{GDA} scheme~\eqref{eqn: implicit-gda}. Let $z = (x, y)^{\top}$ denote the primal-dual variable we define the operator $F(z)$ as
\begin{equation}
\label{eqn: grad-minimax}
F(z) = F(x, y) = \left(\nabla L_{ x}(x, y),  -  \nabla L_{ y}(x, y) \right)^{\top}.
\end{equation}
Consequently, the implicit~\texttt{GDA} scheme~\eqref{eqn: implicit-gda} can be expressed compactly as
\begin{equation}
\label{eqn: implicit-gda-z}
z_{k} = z_{k-1} - s F(z_{k}). 
\end{equation}
Furthermore, the step size $s$ need not remain constant through the iterations, as observed in the proximal~\texttt{ALM}~\eqref{eqn: prox-aug-lag-method}. Let $\{s_{k}\}_{k=1}^{\infty}$ be a variable step-size sequence; for instance, $s_{k} = \gamma_{k} s$, where $\gamma_{k}$ is a polynomial of $k$. Under this variable step-size regime, the implicit \texttt{GDA} scheme~\eqref{eqn: implicit-gda} generalizes naturally to: 
\begin{equation}
\label{eqn: implicit-gda-z-vary}
z_{k} = z_{k-1} - s_{k} F(z_{k}). 
\end{equation}
In this study, this generalization serves as a key observation for deriving the rate $O(1/k)$ for the last iterate. Moreover, we combine this property with second-order ODEs to design new algorithms that further accelerate convergence. Thus, for the corresponding accelerated proximal~\texttt{ALM}s,  their outer iterations successfully achieve these fast convergence rates. 



\subsection{Lyapunov analysis: from objective gap to variational inequality}
\label{subsec: variational-inequality}

In classical convex minimization, Lyapunov analysis typically employs the objective gap, $f(x_k) - f(x^*)$, as a potential function. This choice is natural because the objective gap is always non-negative for convex problems (see, e.g.,~\citet{shi2022understanding}). For minimax problems, however, let  $z^* = (x^*,  y^*)^{\top}$ be the saddle point. In this setting, the objective gap $L(z_k) - L(z^*)$ is generally not sign-definite and therefore cannot serve directly as a Lyapunov measure.  Instead, it is common to utilize the primal-dual gap $L(x_k,  y^*) - L(x^*, y_{k})$. To motivate this viewpoint, we first consider the variational characterization of convex functions. For a convex function $f$, the objective gap satisfies:
\begin{equation}
\label{eqn: variational-inequality-min}
f(x_k) - f(x^*) \leq \langle \nabla f(x_{k}), x_k - x^* \rangle,
\end{equation}
which is a direct consequence of convexity. This variational inequality naturally extends to the minimax setting: 
\begin{equation}
\label{eqn: variational-inequality-minimax}
L(x_k,  y^*) - L(x^*, y_{k}) \leq \left\langle F(z_k), z_{k} - z^* \right\rangle,
\end{equation}
where the monotone operator $F$ is defined in~\eqref{eqn: grad-minimax}.  Motivated by this observation, we replace the conventional objective-gap term in the Lyapunov function with the variational quantity $\left\langle F(z_k), z_{k} - z^* \right\rangle$.  This quantity serves as a natural generalization of the objective gap to saddle-point problems, providing a more suitable potential function for convergence analysis. 

The next step is to investigate the evolution of this quantity along the iterations. In the continuous setting, the time derivative of the objective gap in standard minimization is given by: 
\[
\frac{d \left( f(X) - f(x^*) \right) }{dt} = \langle \nabla f(X), \dot{X} \rangle.
\] 
Analogously, for the minimax problem, the derivative of the variational quantity is:
 \[
 \frac{d}{dt}  \langle F(Z), Z-z^* \rangle =  \langle F(Z), \dot{Z} \rangle + \langle \nabla F(Z) \dot{Z}, Z - z^* \rangle.
 \]  
Compared to the minimization case, the term $ \langle \nabla F(Z) \dot{Z}, Z - z^* \rangle$ arises and requires careful consideration.  A similar phenomenon arises in the discrete setting. While the iterative difference of the objective gap, $f(x_{k+1}) - f(x_k)$, is typically bounded by the variational inequality: 
\begin{equation}
\label{eqn: variationa-inequality}
f(x_{k+1}) - f(x_{k}) \leq  \left\langle \nabla f(x_{k+1}), x_{k+1} - x_{k} \right\rangle, 
\end{equation}
we evaluate the iterative difference of the variational quantity as follows: 
\begin{align}
\left\langle F(z_{k+1}), z_{k+1} - z^* \right\rangle & -  \left\langle F(z_k), z_{k} - z^* \right\rangle \nonumber \\ & =  \left\langle F(z_{k+1}), z_{k+1} - z_{k} \right\rangle +  \left\langle F(z_{k+1}) - F(z_k), z_{k} - z^* \right\rangle \label{eqn: variational-inequality}
\end{align}
The first term is analogous to the objective-gap difference in standard convex optimization. The second term, $ \left\langle F(z_{k+1}) - F(z_k), z_{k} - z^* \right\rangle$, has no counterpart in minimization and represents the core challenge in saddle-point problems. Controlling this residual term is therefore the central component of Lyapunov analysis for minimax algorithms. 

\subsection{Overview of contributions}
\label{subsec: contribution}

In~\Cref{subsec: gda}, we demonstrate the equivalence between the outer iteration of the Proximal~\texttt{ALM}~\eqref{eqn: prox-aug-lag-method} and the implicit~\texttt{GDA} framework~\eqref{eqn: implicit-gda}, including its variable step size extension. In~\Cref{subsec: variational-inequality}, we point out that, in contrast to minimization methods, the focus of the Lyapunov analysis for the \texttt{GDA}-type methods shifts from the objective gap to the variational inequality.  Motivated by these observations, we propose several accelerated implicit \texttt{GDA} schemes. Using continuous-time ODEs and Lyapunov analysis, we prove that the proposed methods achieve faster convergence rates. The main contributions of this study are summarized below. 

\paragraph{An implicit \texttt{GDA} scheme with variable step sizes}  Motivated by a continuous-time ODE perspective, we propose the following implicit GDA scheme with variable step sizes as: 
\begin{equation}
\label{eqn: im-gda-discrete}
z_{k} + (k +2)s  F(z_{k})  = z_{k-1}  + k s F(z_{k-1}),
\end{equation}
where $s>0$ is a fixed parameter. By means of Lyapunov analysis, we establish that this method achieves an $o(1/k)$ convergence rate for the last iterate:
\[
\left\{ \begin{aligned}
        & L(x_k; y^{*}) - L(x^{*}; y_k) \leq o\left( \frac{1}{k} \right), \\
        & \| F(z_k) \| \leq o\left( \frac{1}{k} \right).
        \end{aligned} \right.
\]
This result is derived directly from first-order dynamical system, offering an analytical framework substantially simpler than the second-order dynamics approach proposed by~\citet{boct2025fast}. Furthermore, by setting the step size as $s_k = (k+2)s$, the proximal~\texttt{ALM}~\eqref{eqn: prox-aug-lag-method}, or equivalently the implicit~\texttt{GDA} scheme~\eqref{eqn: implicit-gda-z-vary}, facilitates efficient computation of $z_k$ from  $z_{k-1} + k s F(z_{k-1})$ in practical applications.



\paragraph{A family of Nesterov-type accelerated implicit \texttt{GDA} schemes} Motivated by a second-order ODE perspective, we propose the following family of Nesterov-type accelerated implicit GDA schemes characterized by variable step sizes as: 
\begin{subequations}
 \label{eqn: im-n-gda}
 \begin{empheq}[left=\empheqlbrace]{align} 
         & z_{k} = \xi_{k-1} - \frac{k^r s }{r+2} F(z_k),                                                                                                            \label{eqn: im-n-grad}  \\
         & \xi_{k} = z_{k}  + \frac{k}{k+r+3}\left( z_{k} - z_{k-1} \right) -  \frac{(k+1)^{r+1}}{k^r(k+r+3)} (z_k - \xi_{k-1}),       \label{eqn: im-n-moment}
  \end{empheq}  
 \end{subequations}
 where $s>0$ is a fixed parameter and $r\geq 0$ serves as a family parameter that determines the particular accelerated scheme. Using Lyapunov analysis, we demonstrate that for any given $r$, there exists a threshold $k_0 := k_0(r) > 0$ such that, for all $k \geq k_0$, the proposed scheme achieves the following last-iterate convergence rate:
  \[
L(x_k; y^{*}) - L(x^{*}; y_k) \leq o\left( \frac{1}{k^{r+1}} \right),
\]
which is significantly faster than the $o(1/k)$ last-iterate convergence rate established by~\citet{boct2025fast}. Similarly, by setting the step size as $s_k = sk^r / (r+2)$, the proximal~\texttt{ALM}~\eqref{eqn: prox-aug-lag-method}, or equivalently the implicit~\texttt{GDA} scheme~\eqref{eqn: implicit-gda-z-vary}, facilitates efficient computation of $z_k$ from  $\xi_{k-1} $ in practice.

\paragraph{An explicit \texttt{GDA} scheme with $o(1/k)$ last--iterate rate}
Due to numerical errors inherent in the discretization from the ODE to explicit algorithms, the acceleration parameter must be chosen as $r=0$. Motivated by the resulting second-order dynamical system, we propose the following explicit \texttt{GDA} scheme as:
\begin{subequations}
 \label{eqn: explicit-n-gda}
\begin{empheq}[left=\empheqlbrace]{align} 
 &     \xi_{k} = z_{k-1} - \frac{s}{2} F(z_{k-1})                                                                              \label{eqn: explicit-grad}                 \\
 &     z_{k}   = z_{k-1}   - \frac{s}{k+2} F(z_{k-1})    + \frac{k-1}{k+2} (\xi_{k} - \xi_{k-1})            \label{eqn: explicit-moment}
  \end{empheq}  
 \end{subequations}
which achieves an $o(1/k)$ last--iterate rate for the primal-dual objective gap. Compared with the analysis of~\citet{boct2025fast}, our proof is significantly more transparent and concise.


%
%
%
%


\subsection{Preliminaries and organization}
\label{subsec: organization}

Let $\mathcal{C}^{p}$ denote the space of $p$-times continuously differentiable functions on $\mathbb{R}^n$ for $p\in \{1,2\}$. We define $\mathcal{F}^p \subset \mathcal{C}^{p}$ as the subclass of convex-concave functions possessing a unique saddle point $z^* = (x^*, y^*)^{\top}$ that satisfies:
\begin{equation}
\label{eqn: saddle-point}
L(x^*, y) \leq L(x^*, y^*) \leq L(x, y^*), 
\end{equation}
where the left equality holds if and only if $y = y^*$, and the right equality holds if and only if $x = x^*$. From the property of the monotone operator defined in~\eqref{eqn: grad-minimax}, we derive:
\begin{equation}
\label{eqn: vi-saddle-point-1}
\left\langle F(z), z - z^* \right\rangle \geq L(x, y^*) - L(x^*, y) 
\end{equation}
where the equality holds if and only if $z = z^*$.  Given the convexity of $L$ with respect to $x$ and the concavity of $L$ with respect to $y$, the following inequalities hold for any $z_1, z_2$:
\[
\left\{ \begin{aligned}
         & L(z_2) - L(x_1, y_2)  \leq \left\langle \nabla_x L(z_2), x_2 - x_1 \right\rangle \\
         & L(z_1) - L(x_2, y_1)  \leq \left\langle \nabla_x L(z_1), x_1 - x_2 \right\rangle  
         \end{aligned} \right.
\]
and
\[
\left\{ \begin{aligned}
         & L(x_1, y_2) - L(z_1)   \leq \left\langle \nabla_y L(z_1), y_2 - y_1 \right\rangle  \\
         &  L(x_2, y_1) - L(z_2)  \leq  \left\langle \nabla_y L(z_2), y_1 - y_2 \right\rangle 
         \end{aligned} \right.
\]
Summing these four inequalities, we derive the following variational inequality: 
\begin{equation}
\label{eqn: vi-saddle-point-2}
\left\langle F(z_1) - F(z_2), z_1 - z_2 \right\rangle \geq 0. 
\end{equation}
Furthermore, it is straightforward to show that the Jacobian $\nabla F$ is a matrix with a positive semi-definite symmetric part, which implies that $\nabla F$ is positive semi-definite in the sense that $v^T\nabla F(z)v\geq 0$ for all $v$. According to the uniqueness of the saddle point, we know that $\nabla F(z^*)$ is positive definite. 

The remainder of this paper is organized as follows.~\Cref{sec: implicit-gda} analyzes the implicit \texttt{GDA} scheme with variable step sizes~\eqref{eqn: im-gda-discrete}, while~\Cref{sec: nag} examines the family of Nesterov-type accelerated implicit \texttt{GDA} schemes~\eqref{eqn: im-n-gda}. In both sections, we first derive the corresponding continuous-time ODEs and establish Lyapunov functions to determine convergence rates, subsequently extending these results to the discrete algorithms.~\Cref{sec: explicit} utilizes these continuous-time results to construct a Nesterov-type explicit \texttt{GDA} scheme~\eqref{eqn: explicit-n-gda} and derives the $o(1/k)$ last-iterate convergence rate for the primal-dual objective gap.~\Cref{sec: numerical} shows the numerical experiments for the benchmark test and the steady-state Stokes problem.~\Cref{sec: conclusion} offers concluding remarks and discusses potential directions for future research.



%
%
%
%
%
%
%

\section{An implicit GDA scheme with variable step sizes}
\label{sec: implicit-gda}

In this section, we establish the convergence rates for both the primal-dual objective gap and the gradient norm regarding the implicit \texttt{GDA} scheme with variable step sizes, as given in~\eqref{eqn: im-gda-discrete}. We begin by deriving the corresponding continuous-time ODE and constructing an appropriate Lyapunov function. This enables us to characterize the convergence behavior of the continuous-time dynamics and establish the associated convergence rates. We then extend these results to the discrete-time algorithm~\eqref{eqn: im-gda-discrete}, thereby obtaining rigorous convergence guarantees.


\subsection{Continuous-time dynamics}
\label{subsec: continuous-gda}

We first reformulate the implicit \texttt{GDA} scheme with variable step sizes~\eqref{eqn: im-gda-discrete} as: 
\begin{equation}
\label{eqn: im-gda-discrete-1}
z_{k} - z_{k-1} + ks \left( F(z_{k}) - F(z_{k-1}) \right) = -2s F(z_{k}). 
\end{equation}
Taking the limit as the parameter $s \rightarrow 0$, the finite difference $(z_k - z_{k-1}) / s$ approximates to the derivative $\dot{Z}$. Consequently, we derive the following ODE corresponding to the discrete algorithm~\eqref{eqn: im-gda-discrete-1} as 
\begin{equation}
\label{eqn: im-gda-con}
\dot{Z} + t \nabla F(Z) \dot{Z} = - 2F(Z),
\end{equation}
with any initial $Z(0) = z_0 \in \mathbb{R}^n$. To analyze the convergence. we define the Lyapunov function as: 
\begin{align}
\mathcal{E}(t) & = \frac{1}{2} \| Z + t F(Z) - z^* \|^2 \nonumber \\
                       & = \frac{1}{2} \| Z - z^* \|^2 + t\left\langle F(Z), Z - z^* \right\rangle + \frac{t^2}{2} \| F(z) \|^2. \label{eqn: im-gda-con-ly}
\end{align}
From the variational inequality~\eqref{eqn: vi-saddle-point-1}, we know that $\left\langle F(Z), Z - z^* \right\rangle \geq 0$, which implies $\mathcal{E}(t) \geq 0$. Taking the time derivative, we obtain:
\begin{align}
\frac{d \mathcal{E}(t)}{ dt} & = \big\langle \dot{Z} + t \nabla F(Z) \dot{Z} + F(Z), Z + t F(Z) - z^* \big \rangle           \nonumber                                   \\
                                          & \leq - \big \langle  F(Z), Z  - z^* \big \rangle - t\| F(Z) \|^2  \leq 0,                                  \label{eqn: im-gda-con-ly-der}
\end{align}
where the equality holds if and only if $Z = z^*$, implying that $\mathcal{E}(t) \rightarrow (1/2)\|z^* - tF(z^*) - z^* \|^2 = 0$ as $t \rightarrow \infty$. Thus, we establish the following convergence rates. 
\begin{theorem}
\label{thm: im-gad-con-rate}
For any convex-concave function $L\in \mathcal{F}^2$, the solution $Z = (X, Y)$ to the continuous-time system~\eqref{eqn: im-gda-con} satisfies: 
\begin{subequations}
 \label{eqn: im-gda-con-rate}
        \begin{empheq}[left=\empheqlbrace]{align} 
        & L(X, y^*) - L(x^*, Y) \leq \frac{\| z_0 - z^* \|^2}{2t},        \label{eqn: im-gda-con-rate-obj-gap}              \\
        & \| F(Z) \|^2 \leq \frac{\| z_0 - z^* \|^2}{t^2}.                               \label{eqn: im-gda-con-rate-grad-norm} 
        \end{empheq}
\end{subequations}
Moreover, the following asymptotic properties hold:
\begin{subequations}
 \label{eqn: im-gda-con-rate-o}
        \begin{empheq}[left=\empheqlbrace]{align} 
        & \lim_{t \rightarrow \infty} t \left( L(X, y^*) - L(x^*, Y) \right) = 0,        \label{eqn: im-gda-con-rate-obj-gap-o}              \\
        & \lim_{t \rightarrow \infty} t^2 \| F(Z) \|^2 = 0.                                     \label{eqn: im-gda-con-rate-grad-norm-o} 
        \end{empheq}
\end{subequations}
\end{theorem}

\subsection{Convergence analysis for the discrete algorithm}
\label{subsec: discrete-gda}

To analyze the convergence of the discrete algorithm, we introduce a discrete-time Lyapunov function corresponding to its continuous-time counterpart~\eqref{eqn: im-gda-con-ly-der} as: 
\begin{align}
\mathcal{E}(k) & = \frac{1}{2} \| z_k + ks F(z_k) - z^* \|^2 \nonumber \\
                       & = \frac{1}{2} \| z_k - z^* \|^2 + ks \left\langle F(z_k), z_k - z^* \right\rangle + \frac{k^2s^2}{2} \| F(z_k) \|^2 \label{eqn: im-gda-discrete-ly}
\end{align}
Using the implicit~\texttt{GDA} method with variable step sizes~\eqref{eqn: im-gda-discrete},  the iterative difference in the Lyapunov function satisfies
\begin{equation}
\label{eqn: im-gda-discrete-ly-diff}
\mathcal{E}(k+1) -  \mathcal{E}(k) \leq - s \big\langle F(z_{k+1}), z_{k+1} - z^* \big \rangle - \frac{(2k+3)s^2}{2} \| F(z_{k+1}) \|^2, 
\end{equation}
where the equality holds if and only if $z_{k+1} = z^*$. Consequently, as the iteration progress, $k \rightarrow \infty$, the Lyapunov function satisfies $\mathcal{E}(k) \rightarrow (1/2)\|z^* - tF(z^*) - z^* \|^2 = 0$. Based on this descent property, we establish the following convergence rates.
\begin{theorem}
\label{thm: im-gad-dis-rate}
For any convex-concave function $L\in \mathcal{F}^1$, the iterative sequence $\{ z_k\}_{k=0}^{\infty}  = \{ (x_k, y_k) \}_{k=0}^{\infty}$ generated by the implicit~\texttt{GDA} method with variable step sizes~\eqref{eqn: im-gda-discrete} satisfies:
\begin{subequations}
 \label{eqn: im-gda-dis-rate}
        \begin{empheq}[left=\empheqlbrace]{align} 
        & L(x_k, y^*) - L(x^*, y_k) \leq \frac{\| z_0 - z^* \|^2}{2ks},                    \label{eqn: im-gda-dis-rate-obj-gap}              \\
        & \| F(z_k) \|^2 \leq \frac{\| z_0 - z^* \|^2}{k^2s^2}.                                \label{eqn: im-gda-dis-rate-grad-norm} 
        \end{empheq}
\end{subequations}
Furthermore, the following asymptotic properties hold:
\begin{subequations}
 \label{eqn: im-gda-dis-rate-o}
        \begin{empheq}[left=\empheqlbrace]{align} 
        & \lim_{k \rightarrow \infty} k \left( L(x_k, y^*) - L(x^*, y_k) \right) = 0,        \label{eqn: im-gda-dis-rate-obj-gap-o}              \\
        & \lim_{k \rightarrow \infty} k^2 \| F(z_k) \|^2 = 0.                                        \label{eqn: im-gda-dis-rate-grad-norm-o} 
        \end{empheq}
\end{subequations}
\end{theorem}


\section{A family of Nesterov-type accelerated implicit GDA schemes}
\label{sec: nag}

In this section, we establish the convergence rates for the primal-dual objective gap within the family of Nesterov-type accelerated implicit~\texttt{GDA} schemes defined  in~\eqref{eqn: im-n-gda}. We begin by deriving the corresponding family of second-order ODEs and constructing appropriate Lyapunov functions, which allows us to characterize the convergence behavior of the continuous-time dynamics. Building upon these insights, we extend our analysis to the discrete-time algorithm~\eqref{eqn: im-n-gda} to establish rigorous convergence guarantees.

\subsection{Continuous-time dynamics for acceleration}
\label{subsec: continuous-nag}

We begin by substituting the gradient iteration~\eqref{eqn: im-n-grad} into the momentum iteration~\eqref{eqn: im-n-moment} to derive the relationship for the iterative sequence $\{ z_{k} \}_{k=0}^{\infty}$ as:
\begin{align}
& (k+1)(z_{k+1} - z_{k})  +  (r+2)(z_{k+1} - z^*)  \nonumber \\ 
& = k(z_{k} - z_{k-1}) + (r+2)(z_{k} - z^*) -  (k+1)^{r} s F(z_{k+1}) - \frac{  \left(k+1 \right)^{r+1} s}{r+2} \left( F(z_{k+1}) - F(z_{k}) \right).   \label{eqn: implicit-nag}
\end{align}
Taking the step size $\rho = s^{\frac{1}{r+1}}$,we define the time variable as $t=k\rho$. At the step size approaches zero ($\rho \rightarrow 0$), we derive the following second-order ODE as
\begin{equation}
\label{eqn: im-n-ode}
\ddot{Z} + \frac{r+3}{t} \dot{Z} + \frac{  t^{r}}{r+2} \nabla F(Z) \dot{Z} +  t^{r-1}  F(Z) = 0,
\end{equation}
with any initial $Z(0) = z_0 \in \mathbb{R}^n$ and $\dot{Z}(0) = 0$. To analyze the convergence. we propose the following Lyapunov function as: 
\begin{equation}
\label{eqn: im-n-ode-lyapunov}
\mathcal{E}(t) = \underbrace{ t^{r+1} \left\langle F(Z), Z - z^{*} \right\rangle }_{:= \; \mathcal{E}_1(t)} + \underbrace{ \frac{1}{2} \big \| t \dot{Z} + (r+2)(Z - z^*) \big\|^2}_{: = \; \mathcal{E}_2(t)}.
\end{equation}
Taking the time derivatives for the potential energy $\mathcal{E}_1(t)$ and the mixed energy $\mathcal{E}_2(t)$, we obtain:
\begin{equation}
\label{eqn: im-n-ode-lyapunov-der-1}
\frac{d \mathcal{E}_1(t)}{dt}  
         = \underbrace{(r + 1) t^{r} \big\langle F(Z), Z - x^{*} \big\rangle}_{\mathbf{I}_1} + \underbrace{ t^{r + 1} \big \langle \nabla F(Z) \dot{Z}, Z - z^{*} \big \rangle}_{\mathbf{I}_2} + \underbrace{ t^{r + 1} \big \langle F(Z),  \dot{Z} \big \rangle}_{\mathbf{I}_3}
\end{equation}
and
\begin{align}
\frac{d \mathcal{E}_2(t)}{dt} & =  \big\langle t\ddot{Z} + (r + 3) \dot{Z},  t \dot{Z} + (r+2)(Z - z^{*}) \big\rangle                                                   \nonumber                \\
                                             & = -  \left\langle  \frac{t^{r+1}}{r+2} \nabla F(Z) \dot{Z} + t^{r} F(Z),   t \dot{Z} + (r+2)(Z - z^{*})  \right\rangle      \nonumber                \\
                                             & = - \underbrace{(r + 2)  t^{r} \big\langle F(Z), Z - z^{*} \big\rangle}_{\mathbf{II}_1} - \underbrace{ t^{r + 1} \big \langle \nabla F(Z) \dot{Z}, Z - z^{*} \big \rangle}_{\mathbf{II}_2} - \underbrace{ t^{r + 1} \big \langle F(Z),  \dot{Z} \big \rangle}_{\mathbf{II}_3}                             \nonumber                \\
                                             & \mathrel{\phantom{=}}  - \frac{ t^{r+2}}{r + 2} \big\langle \nabla F(Z) \dot{Z}, \dot{Z} \big\rangle.                        \label{eqn: im-n-ode-lyapunov-der-2}
\end{align}
Comparing~\eqref{eqn: im-n-ode-lyapunov-der-1} and~\eqref{eqn: im-n-ode-lyapunov-der-2}, we observe that $\mathbf{I}_2 = \mathbf{II}_2$, $\mathbf{I}_3 = \mathbf{II}_3$, as well as $\mathbf{I}_1 - \mathbf{II}_1 = - t^{r} \big\langle F(Z), Z - x^{*} \big\rangle$. Thus, the derivative of the Lyapunov function satisfies: 
\begin{equation}
\label{eqn: im-n-ode-lyapunov-der-3}
\frac{d \mathcal{E}(t)}{dt}  =  - t^{r + 1} \big \langle F(Z),  Z - z^* \big \rangle  - \frac{ t^{r+2}}{r + 2} \big\langle \nabla F(Z) \dot{Z}, \dot{Z} \big\rangle \leq 0. 
\end{equation}
Since $\nabla F(z^*)$ is positive definite, the above equality holds if and only if $Z = z^*$ and $\dot{Z} = 0$, implying that $\mathcal{E}(t) \rightarrow \mathcal{E}(Z = z^*) = 0$ as $t \rightarrow \infty$. Consequently, we establish the following convergence rates. 
\begin{theorem}
\label{thm: im-n-ode-rate}
For any convex-concave function $L\in \mathcal{F}^2$, the solution $Z = (X, Y)$ to the continuous-time system~\eqref{eqn: im-n-ode} satisfies: 
\begin{equation}
 \label{eqn: im-n-ode-rate}
 L(X, y^*) - L(x^*, Y) \leq \frac{ (r+2)^2\| z_0 - z^* \|^2}{2t^{r+1}}.              
\end{equation}
Moreover, the following asymptotic properties hold:
\begin{equation}
 \label{eqn: im-n-ode-rate-o}
 \lim_{t \rightarrow \infty} t^{r+1} \left( L(X, y^*) - L(x^*, Y) \right) = 0.     
\end{equation}
\end{theorem}

\subsection{Accelerated convergence for the family of discrete algorithms}
\label{subsec: discrete-nag}

To analyze the convergence of the discrete algorithm, we introduce a discrete-time Lyapunov function corresponding to its continuous-time counterpart~\eqref{eqn: im-n-ode-lyapunov} as: 
\begin{equation}
\label{eqn: im-n-dis-lyapunov}
\mathcal{E}(k) = \underbrace{\left(k+1\right)^{r+1}s \big\langle F(z_k), z_k - z^{*} \big\rangle}_{: = \; \mathcal{E}_1(k)} + \underbrace{ \frac{1}{2} \big \| k(z_{k} - z_{k-1}) + (r+2)(z_k - z^{*}) \big \|^2 }_{: =  \; \mathcal{E}_2(k)}.
\end{equation}
Calculating the iterative difference for the potential energy $\mathcal{E}_1(k)$, we obtain:
\begin{align}
\mathcal{E}_1(k+1) - \mathcal{E}_1(k) = & \underbrace{ \left[ (k+2)^{r+1} - (k+1)^{r+1} \right] s  \big\langle F(z_{k+1}), z_{k+1} - z^{\star} \big\rangle}_{\mathbf{I}_1}   \nonumber \\ 
                                                                 & + \underbrace{ \left(k+1 \right)^{r+1} s  \big\langle F(z_{k+1}) - F(z_{k}), z_{k+1} - z^{\star} \big\rangle } _{\mathbf{I}_2}      \nonumber \\
                                                                 & + \underbrace{  \left(k+1 \right)^{r+1}  s \big\langle F(z_{k}), z_{k+1} - z_{k} \big\rangle }_{\mathbf{I}_3 },                              \label{eqn: im-n-dis-lyapunov-der-1}
\end{align}
where the coefficient of the term $\mathbf{I}_1$ satisfies the estimate as $(k+2)^{r+1} - (k+1)^{r+1} \leq (r+1)(k+1)^r + 2^{r+1}(k+1)^{r-1}$. Furthermore, utilizing the the relationship~\eqref{eqn: implicit-nag}, we calculate the iterative difference for the mixed energy $\mathcal{E}_2(k)$ as: 
\begin{align}
& \mathcal{E}_2(k+1)  - \mathcal{E}_2(k)                                                                                                                                                                                      \nonumber      \\
& =  - (k+1)^{r} s \left\langle  F(z_{k+1}) + \frac{  k+1 }{r+2} \left(F(z_{k+1}) - F(z_{k}) \right),   (k+1)(z_{k+1} - z_{k})  +  (r+2)(z_{k+1} - z^*)  \right\rangle         \nonumber      \\
& \mathrel{\phantom{=}} - \frac{(k+1)^{2r} s^2 }{2} \left\| F(z_{k+1}) + \frac{  k+1 }{r+2} \left( F(z_{k+1}) - F(z_{k}) \right)  \right\|^2                                             \nonumber      \\
& \leq - \underbrace{ (r+2)(k+1)^{r} s \big\langle F(z_{k+1}), z_{k+1} - z^{*} \big\rangle } _{\mathbf{II}_1}                                                                                    \nonumber      \\
& \mathrel{\phantom{\leq}}  - \underbrace{ \left(k+1 \right)^{r+1} s  \big\langle F(z_{k+1}) - F(z_{k}), z_{k+1} - z^{*} \big\rangle } _{\mathbf{II}_2}                      \nonumber       \\
& \mathrel{\phantom{\leq}}  - \underbrace{  \left(k+1 \right)^{r+1} s \big\langle F(z_{k+1}), z_{k+1} - z_{k} \big\rangle }_{\mathbf{II}_3 }                                      \nonumber      \\
&  \mathrel{\phantom{\leq}} - \frac{(k+1)^{r+2} s}{r+2}  \big\langle F(z_{k+1}) - F(z_{k}), z_{k+1} - z_{k} \big\rangle.                                                                     \label{eqn: im-n-dis-lyapunov-der-2}
\end{align}
Comparing~\eqref{eqn: im-n-dis-lyapunov-der-1} and~\eqref{eqn: im-n-dis-lyapunov-der-2}, we observe that $\mathbf{I}_2 = \mathbf{II}_2$ and $\mathbf{I}_3 = \mathbf{II}_3$. Furthermore, the difference between $\mathbf{I}_1$ and $\mathbf{II}_1$ satisfies the inequality $\mathbf{I}_1 - \mathbf{II}_1 \leq - (k+1 - 2^{r+1}) (k+1)^{r-1} s \big\langle F(z_k), z_k - x^{*} \big\rangle$. Consequently, the iterative difference of the Lyapunov function satisfies: 
\begin{align}
\mathcal{E}(k+1) - \mathcal{E}(k) \leq &  - (k+1 - 2^{r+1})(k+1)^{r-1} s  \big\langle F(z_{k+1}), z_{k+1} - z^{*} \big\rangle   \nonumber \\ &  - \frac{(k+1)^{r+2} s}{r+2}  \big\langle F(z_{k+1}) - F(z_{k}), z_{k+1} - z_{k} \big\rangle \leq 0 \label{eqn: im-n-dis-lyapunov-der-3}
\end{align}
From the variational inequality~\eqref{eqn: vi-saddle-point-1}, we know that $ \big\langle F(z_{k+1}), z_{k+1} - z^{*} \big\rangle =0$ holds if and only if $z_{k+1} = z^*$.  Given the structure of the difference~\eqref{eqn: im-n-dis-lyapunov-der-3}, $z_{k+1} = z^*$ implies $z_k = z^*$, which ensures that $\mathcal{E}(k) \rightarrow \mathcal{E}(z_k = z^*) = 0$ as $k \rightarrow \infty$. This monotonicity and convergence to zero allow us to establish the following convergence rates.
\begin{theorem}
\label{thm: im-n-dis-rate}
For any convex-concave function $L\in \mathcal{F}^1$, let $\{ z_k\}_{k=0}^{\infty}  = \{ (x_k, y_k) \}_{k=0}^{\infty}$ be the iterative sequence generated by the family of Nesterov-type accelerated implicit~\texttt{GDA} schemes~\eqref{eqn: im-n-gda}. For any given $r \geq 0$, there exists a threshold $k_0 := k_0(r) = 2^{r +1}> 0$ such that, for all $k \geq k_0$, the following bound holds:
\begin{equation}
 \label{eqn: im-n-dis-rate}
 L(x_k, y^*) - L(x^*, y_k) \leq \frac{\mathcal{E}(k_0)}{2(k+1)^{r+1}s}.             
\end{equation}
Furthermore, the sequence satisfies the following asymptotic property:
\begin{equation}
 \label{eqn: im-n-dis-rate-o}
 \lim_{k \rightarrow \infty} (k+1)^{r+1} \left( L(x_k, y^*) - L(x^*, y_k) \right) = 0.
 \end{equation}
\end{theorem}

\section{An explicit GDA scheme with $o(1/k)$ last--iterate  rate}
\label{sec: explicit}

%

%
%


In this section, we establish the convergence rate for the primal-dual objective gap under the explicit~\texttt{GDA} scheme defined  in~\eqref{eqn: explicit-n-gda}. By advancing the momentum iteration \eqref{eqn: explicit-moment} by one step and substituting the gradient iteration \eqref{eqn: explicit-grad} into it, we obtain the following recurrence relation:
\begin{align}
(k+1)(z_{k+1} - z_{k}) + 2(z_{k+1} - z^*)  = & k(z_{k} - z_{k-1}) + 2(z_{k} - z^*)  \nonumber \\
                                                                    & -   s F(z_{k}) - \frac{ks}{2} \left( F(z_k) - F(z_{k-1}) \right).   \label{eqn: explicit-nag}
\end{align}
Taking the step size $\rho = s$, we define the continuous-time variable $t=k\rho$,  As the step size approaches zero ($\rho \rightarrow 0$), we derive the following second-order ODE as
\begin{equation}
\label{eqn: im-e-ode}
\ddot{Z} + \frac{3}{t} \dot{Z} + \frac{  1}{2} \nabla F(Z) \dot{Z} +  \frac{1}{t}  F(Z) = 0,
\end{equation}
with any initial $Z(0) = z_0 \in \mathbb{R}^n$ and $\dot{Z}(0) = 0$.   Note that the second-order ODE~\eqref{eqn: im-e-ode} is a special instance of the general family~\eqref{eqn: im-n-ode} where $r=0$. Consequently, the Lyapunov function corresponds to the special case of the family~\eqref{eqn: im-n-ode-lyapunov} with $r=0$, given by: 
\begin{equation}
\label{eqn: im-e-ode-lyapunov}
\mathcal{E}(t) = \underbrace{ t \left\langle F(Z), Z - z^{*} \right\rangle }_{:= \; \mathcal{E}_1(t)} + \underbrace{ \frac{1}{2} \big \| t \dot{Z} + 2(Z - z^*) \big\|^2}_{: = \; \mathcal{E}_2(t)},
\end{equation}
which leads to the convergence rates provided in~\eqref{eqn: im-n-ode-rate} and~\eqref{eqn: im-n-ode-rate-o} for the case where $r=0$. Unlike implicit schemes, we must account for numerical errors to bridge the gap between the continuous-time ODE to its explicit discretization. Although the operator $F$ is monotone by virtue of the variational inequality~\eqref{eqn: vi-saddle-point-2}, we must impose a stricter condition: the inverse operator $F^{-1}$ must satisfy strong monotonicity, as defined by~\citet[Definition 12.53]{rockafellar1998variational}. This assumption is formally stated below.
\begin{defn}
\label{defn: strong-monotonicity}
The operator $F$ is said to have an inverse that satisfies strong monotonicity if the following variational inequality holds:
\begin{equation}
\label{eqn: m-coercive}
 \left\langle F(z_1) - F(z_2), z_1 - z_2  \right\rangle \geq \frac{1}{M} \| F(z_1) - F(z_2) \|^2, 
\end{equation}
for any $z_{1}, z_{2} \in \mathbb{R}^{n} \times \mathbb{R}^{n}$.
\end{defn}
To analyze the convergence of the discrete algorithm, we introduce a discrete-time Lyapunov function corresponding to its continuous-time counterpart~\eqref{eqn: im-e-ode-lyapunov} as: 
\begin{equation}
\label{eqn: im-n-dis-lyapunov}
\mathcal{E}(k) = \underbrace{sk\left\langle F(z_{k-1}), z_{k-1} - z^* \right\rangle }_{:=\;\mathcal{E}_1(k)}  + \underbrace{ \frac{1}{2} \left\| k(z_{k} - z_{k-1}) + 2(z_{k} - z^*)  \right\|^2 }_{:=\;\mathcal{E}_2(k)}.
\end{equation}
For the potential energy $\mathcal{E}_1(k)$, we calculate its iterative difference as 
\begin{align}
\mathcal{E}_1(k+1) - \mathcal{E}_1(k) = & \underbrace{s \left\langle F(z_{k}), z_{k} - z^* \right\rangle}_{\mathbf{I}_1} + \underbrace{sk \left\langle F(z_{k}), z_{k} - z_{k-1} \right\rangle}_{\mathbf{II}_1} \nonumber \\ 
                                                                 & + \underbrace{sk  \left\langle F(z_{k}) - F(z_{k-1}),  z_{k-1} - z^* \right\rangle}_{\mathbf{III}_1}.          \label{eqn: im-e-dis-lyapunov-der-1}
\end{align}
Furthermore, utilizing the the relationship~\eqref{eqn: explicit-nag}, we calculate the iterative difference for the mixed energy $\mathcal{E}_2(k)$ as: 
\begin{align}
 \mathcal{E}_2(k+1) &- \mathcal{E}_2(k) \nonumber \\
& = - s \left\langle  F(z_{k}) +  \frac{k}{2}\left( F(z_k) - F(z_{k-1}) \right), k(z_{k} - z_{k-1}) + 2(z_{k} - z^*)  \right\rangle  \nonumber \\
& \mathrel{\phantom{=}} + \frac{s^2}{2} \left\| F(z_{k}) + \frac{k}{2} \left( F(z_k) - F(z_{k-1}) \right) \right\|^2                    \nonumber \\
&= - \underbrace{ 2 s \left\langle F(z_{k}), z_{k} - z^* \right\rangle}_{\mathbf{I}_2} - \underbrace{sk \left\langle F(z_{k}), z_{k} - z_{k-1} \right\rangle}_{\mathbf{II}_2} - \underbrace{sk  \left\langle F(z_{k}) - F(z_{k-1}),  z_{k} - z^* \right\rangle}_{\mathbf{III}_2} \nonumber \\
& \mathrel{\phantom{=}} - \frac{sk^2}{2}  \left\langle F(z_k) - F(z_{k-1}), z_k - z_{k-1} \right\rangle  + \frac{s^2}{2} \left\| F(z_{k}) + \frac{k}{2} \left( F(z_k) - F(z_{k-1}) \right) \right\|^2.  \label{eqn: im-e-dis-lyapunov-der-2}
\end{align}
By comparing~\eqref{eqn: im-e-dis-lyapunov-der-1} and~\eqref{eqn: im-e-dis-lyapunov-der-2} and invoking the inverse strong monotonicity of $F$, we first derive the following inequalities as 
\begin{subequations}
 \label{eqn: exp-inq-1}
        \begin{empheq}[left=\empheqlbrace]{align} 
        & \mathbf{I}_1 - \mathbf{I}_2 = - s \left\langle F(z_{k}), z_{k} - z^* \right\rangle \leq - \frac{s}{M} \| F(z_k) \|^2,                                                      \label{eqn: exp-inq-1-1}              \\
        & - \frac{sk^2}{2}  \left\langle F(z_k) - F(z_{k-1}), z_k - z_{k-1} \right\rangle \leq -\frac{sk^2}{2M} \| F(z_k) - F(z_{k-1}) \|^2.                                        \label{eqn: exp-inq-1-2} 
        \end{empheq}
\end{subequations}
Noting that $\mathbf{II}_1 = \mathbf{II}_2$,  and considering the remaining terms, we observe that the monotonicity of the operator $F$ implies: 
\begin{equation}
\label{eqn: exp-inq-2}
\mathbf{III}_1 - \mathbf{III}_2 = - sk \left\langle F(z_{k}) - F(z_{k-1}),  z_{k} - z_{k-1} \right\rangle \leq 0. 
\end{equation}
Combining~\eqref{eqn: exp-inq-1-1},~\eqref{eqn: exp-inq-1-2}, and~\eqref{eqn: exp-inq-2}, when the parameter satisfies $s \leq 1/M$, we obtain the iterative difference satisfies:
\begin{align}
\mathcal{E}(k+1) - \mathcal{E}(k) \leq     & - \frac{3s^2k^2}{8}  \left\| F(z_{k}) - F(z_{k-1}) \right\|^2  +  \frac{s^2k}{2} \left\langle  F(z_{k}) - F(z_{k-1}), F(z_{k})  \right\rangle          \nonumber \\
& - \frac{s^2}{2}  \| F(z_k) \|^2 - sk \left\langle F(z_{k}) - F(z_{k-1}),  z_{k} - z_{k-1} \right\rangle                \nonumber \\
                                                     \leq  & - \frac{s^2}{3}  \| F(z_k) \|^2 - sk \left\langle F(z_{k}) - F(z_{k-1}),  z_{k} - z_{k-1} \right\rangle.                                                              \label{eqn: im-e-dis-lyapunov-der-3}
\end{align}
Since $ \| F(z_k) \| = 0$ holds if and only if $z_{k} = z^*$ and given the variational inequality~\eqref{eqn: vi-saddle-point-1}, $z_{k} = z^*$ implies $z_{k-1} = z^*$, which ensures that $\mathcal{E}(k) \rightarrow \mathcal{E}(z_k = z^*) = 0$ as $k \rightarrow \infty$. This monotonicity and convergence to zero allow us to establish the following convergence rates.
\begin{theorem}
\label{thm: ex-n-dis-rate}
For any convex-concave function $L\in \mathcal{F}^1$ satisfying the inverse strong monotonicity (as given in~\Cref{defn: strong-monotonicity}), the iterative sequence $\{ z_k\}_{k=0}^{\infty}  = \{ (x_k, y_k) \}_{k=0}^{\infty}$ generated by the explicit~\texttt{GDA} scheme~\eqref{eqn: explicit-n-gda} satisfies:
\begin{equation}
 \label{eqn: ex-n-dis-rate}
 L(x_k, y^*) - L(x^*, y_k) \leq \frac{2\| z_0 - z^* \|^2}{ks}.             
\end{equation}
Furthermore, the following asymptotic properties hold:
\begin{equation}
 \label{eqn: ex-n-dis-rate-o}
 \lim_{k \rightarrow \infty} k \left( L(x_k, y^*) - L(x^*, y_k) \right) = 0.
 \end{equation}
\end{theorem}

\begin{remark}
For the constrained convex optimization problem~\eqref{eqn: cvx-equal-constraint}, the explicit~\texttt{GDA} scheme~\eqref{eqn: explicit-nag} can be applied directly to solve the augmented Lagrangian function~\eqref{eqn: aug-lagrangian}, thereby avoiding the use of a proximal implicit operator. However, this explicit approach possesses an inherent limitation: it strictly requires additional regularity conditions, such as the Lipschitz gradient condition or the inverse strong monotonicity condition~(\Cref{defn: strong-monotonicity}), to ensure stability and convergence. In practice, the step size must satisfy $s \leq 1/M$. Since the Lipschitz constant $M$ is often unknown, one is forced to select a conservative, very small step size, which inherently leads to slow convergence. It is worth noting that while our analysis utilizes the inverse strong monotonicity condition, this is slightly weaker than the Lipschitz gradient condition employed in the work of~\citet{boct2025fast}, and thus it allows for a significantly more concise and transparent proof of convergence. Alternatively, for the constrained convex optimization problem~\eqref{eqn: cvx-equal-constraint}, one can employ the projected gradient method:
\begin{align}
x_{k} = \mathrm{Proj}_{\{x | Ex = h\}}\left(x_{k-1} - s \nabla f(x_{k-1}) \right)    =  x_{k-1} - s \left(I - E^{\top}(EE^{\top}) E \right) \nabla f(x_{k-1}),
\end{align}
which achieves a last-iterate convergence rate of $O(1/k)$. By incorporating Nesterov-style acceleration, we obtain:
\begin{equation}
\label{eqn: nag}
\left\{ \begin{aligned}
         & y_{k} = \mathrm{Proj}_{\{x | Ex = h\}}\left(x_{k-1} - s \nabla f(x_{k-1}) \right)  \\
         & x_{k} = y_{k} + \frac{k-1}{k+2} (y_{k} - y_{k-1})
        \end{aligned} \right.
\end{equation}
which improves the last-iterate convergence rate to $O(1/k^2)$. Despite these theoretical advantages, determining a suitable step size remains a practical challenge. Furthermore, these methods also require the computation of the matrix inverse $(EE^\top)^{-1}$; although the dimension of the constraint space is often smaller than $n$, this still constitutes a non-negligible computational cost.
\end{remark}

\section{Numerical experiments}
\label{sec: numerical}

In this section, we demonstrate the numerical performance of the implicit~\texttt{GDA} scheme~\eqref{eqn: im-gda-discrete} and the family of Nesterov-type accelerated implicit GDA schemes~\eqref{eqn: im-n-gda}. These schemes are implemented using the outer iterations of their corresponding proximal~\texttt{ALM}s to solve convex optimization problems subject to linear equality constraints~\eqref{eqn: cvx-equal-constraint}. For brevity, we denote these algorithms as~\texttt{PALM} and~\texttt{A-PALM}, where $r$ denotes the parameter, throughout the experiments.


\subsection{A benchmark test}
\label{subsec: qp-linear}

\begin{figure}[htb!]
    \centering
    \begin{subfigure}[t]{0.46\linewidth}
    \centering
    \includegraphics[scale=0.14]{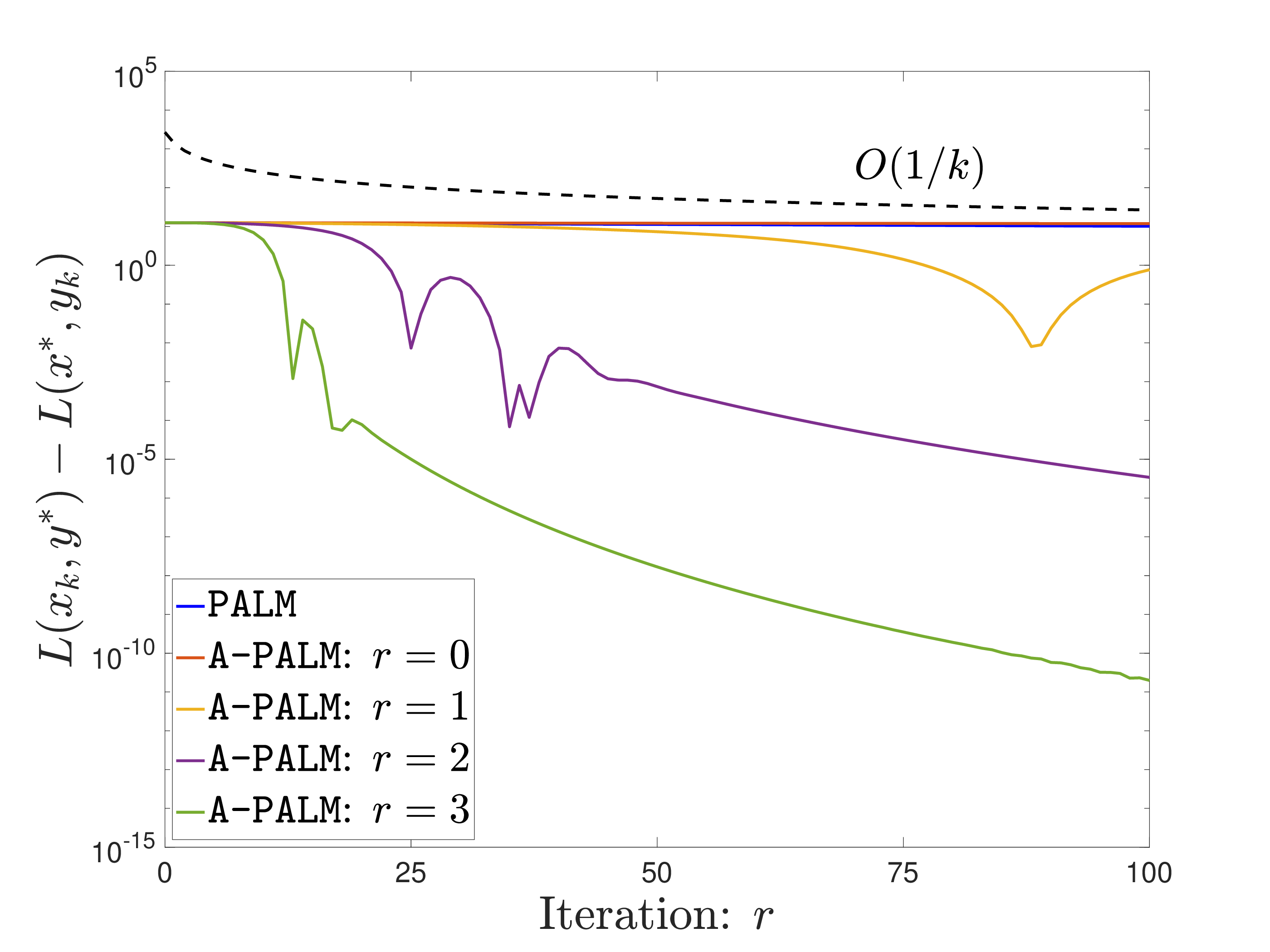}
    \caption{Primal-dual objective gap}
    \label{subfig: obj-benchmark}
    \end{subfigure}
    \begin{subfigure}[t]{0.46\linewidth}
    \centering
    \includegraphics[scale=0.14]{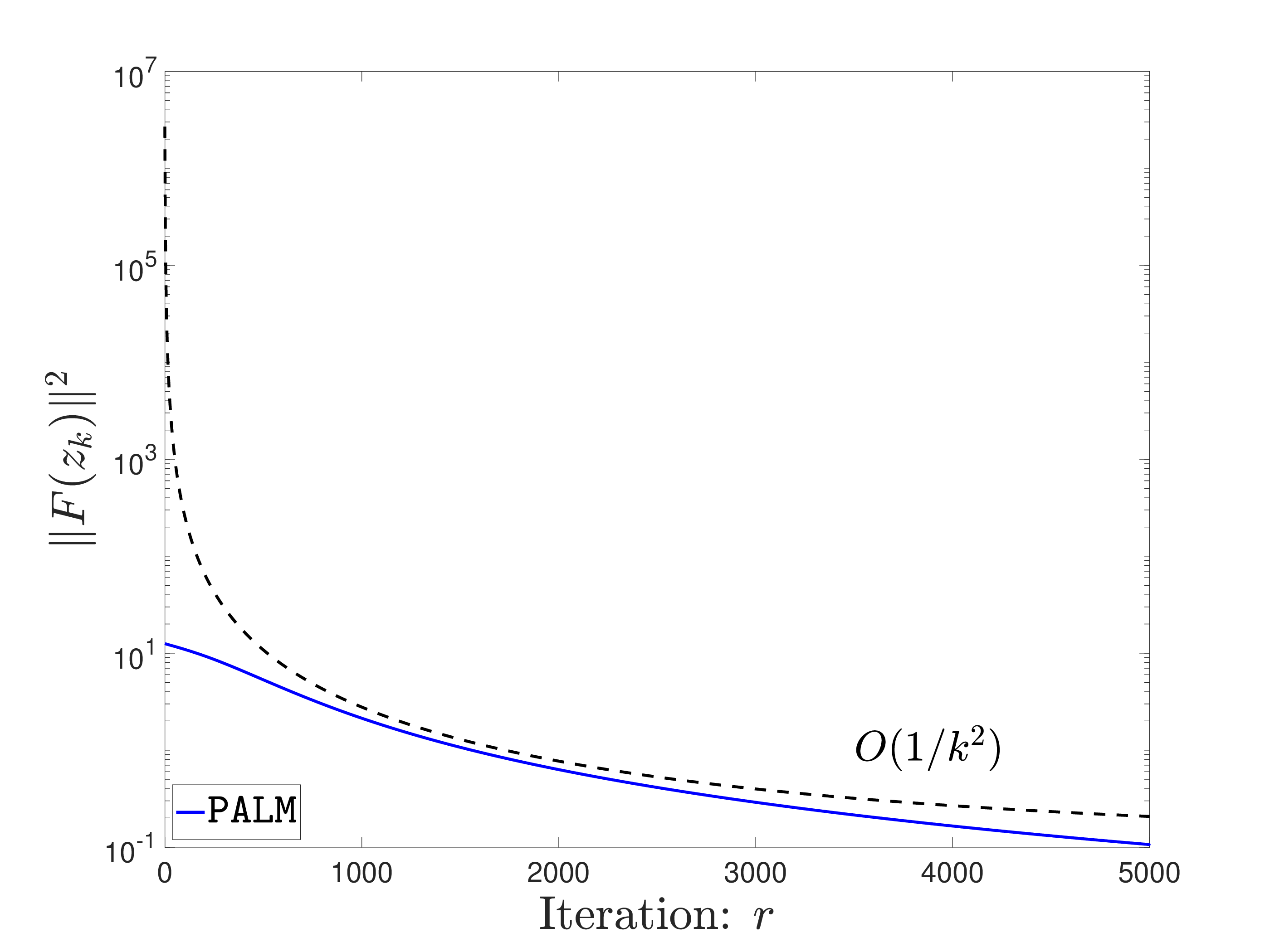}
    \caption{Squared gradient norm}
    \label{subfig: gn-benchmark}
    \end{subfigure}
   \caption{Convergence behavior of~\texttt{PALM} and~\texttt{A-PALM} on the benchmark test. Panel (a) shows the primal-dual objective gap for both algorithms, with~\texttt{A-PALM} set for parameters $r  \in \{0, 1, 2, 3\}$. Panel (b) shows the squared gradient norm for \texttt{PALM}.}
   \label{fig: benckmark}
\end{figure}

To demonstrate the performance of~\texttt{PALM} and~\texttt{A-PALM}, we first conduct numerical experiments based on the setup established in~\citep{ouyang2021lower} and~\citep{yoon2021accelerated}. The benchmark test involves a convex optimization problem subject to the linear equality constraints as defined in~\eqref{eqn: cvx-equal-constraint}. For the least-squares objective function $(1/2) \| Ax - b \|^2$, the matrix $A \in \mathbb{R}^{n \times n}$ and the vector $b \in \mathbb{R}^n$ are defined as:
\[
A = \frac{\sqrt{2}}{4}\begin{pmatrix}
 & & & -1 & 1 \\
 & &  -1& 1 &  \\
 & \iddots& \iddots&  & \\
-1 & 1 & & &\\
1 & & & &
\end{pmatrix}  \quad \mathrm{and} \quad b = \frac{1}{\sqrt{2}} \begin{pmatrix}
1\\ 0 \\  \vdots \\ 0 \\ 0 
\end{pmatrix}.
\]
For the linear equality constraint $Ex = h$, the matrix $E \in \mathbb{R}^{n \times n}$ and the vector $h \in \mathbb{R}^n$ are given as:
\[
E = \frac{1}{4}\begin{pmatrix}
 & & & -1 & 1 \\
 & &  -1& 1 &  \\
 & \iddots& \iddots&  & \\
-1 & 1 & & &\\
1 & & & &
\end{pmatrix} \quad \mathrm{and} \quad  h = \frac{1}{4} \begin{pmatrix}
1\\ 1 \\  \vdots \\ 1 \\ 1 
\end{pmatrix}.
\]
The convergence behavior of~\texttt{PALM} and~\texttt{A-PALM} with parameters $r  \in \{0, 1, 2, 3\}$ is illustrated in~\Cref{fig: benckmark}. 
\Cref{subfig: obj-benchmark} depicts the evolution of the primal-dual objective gap across iterations.  We observe that \texttt{A-PALM} with $r=0$ exhibits convergence nearly identical to that of~\texttt{PALM}, achieving a rate slightly faster than $O(1/k)$. As the parameter $r$ increases,~\texttt{A-PALM} demonstrates substantially accelerated convergence; notably, at $r = 3$, the primal-dual objective gap drops below $10^{-10}$ within only 100 iterations. In~\Cref{subfig: gn-benchmark}, we display the squared gradient norm of~\texttt{PALM} over $5000$ iterations, which confirms the squared gradient norm decays at a rate faster than $O(1/k^2)$ in the long term. Overall, these numerical results are in strong agreement with the theoretical convergence guarantees established in~\Cref{thm: im-gad-dis-rate} and~\Cref{thm: im-n-dis-rate}.



\subsection{Application to the steady-state Stokes system}
\label{subsec: stokes}

We now extend our numerical experiments from the benchmark test to the steady-state Stokes system~\citep{elman2014finite}. Let $\Omega = \{ (x_1, x_2) \in \mathbb{R}^2 | x_1^2 + x_2^2 < 1 \}$ denote the unit disk. Given the body force  $\pmb{f} = (-x_2, x_1)$, the steady-state Stokes system in the disk $\Omega$ is defined as:
\begin{subequations}
 \label{eqn: stokes}
        \begin{empheq}[left=\empheqlbrace]{align} 
         & - \Delta \pmb{u} +  \nabla p = \pmb{f},            \label{eqn: stokes-momentum} \\
         & \nabla \cdot \pmb{u} = 0,                                \label{eqn: stokes-mass} 
        \end{empheq}
\end{subequations}
where $\pmb{u} = (u_1, u_2) \in \mathbb{R}^2$ denotes the velocity field and $p$ represents the pressure. The system is subject to the homogeneous Dirichlet boundary condition 
\begin{equation}
\label{eqn: stokes-boundary}
\pmb{u} = 0, \quad \mathrm{on}\;\; \partial \Omega. 
\end{equation}
By incorporating the Dirichlet boundary condition~\eqref{eqn: stokes-boundary} and utilizing the incompressibility condition~\eqref{eqn: stokes-mass}, we can equivalently transform the steady-state Stokes system~\eqref{eqn: stokes} into the following convex optimization problem with a linear equality constraint:
\begin{subequations}
 \label{eqn: stokes-opt}
        \begin{empheq}[left=\empheqlbrace]{align} 
         & \min J(\pmb{u}) = \int_{\Omega} \| \nabla \pmb{u} \|_2^2 d \pmb{x}- \int_{\Omega} \pmb{f} \cdot \pmb{u} d \pmb{x},    \label{eqn: stokes-obj} \\
         & \mathrm{s.t.} \quad \nabla \cdot \pmb{u} = 0.                                                                                                                 \label{eqn: stokes-con} 
        \end{empheq}
\end{subequations}

\begin{figure}[htb!]
    \centering
     \begin{subfigure}[t]{0.46\linewidth}
    \centering
    \includegraphics[scale=0.14]{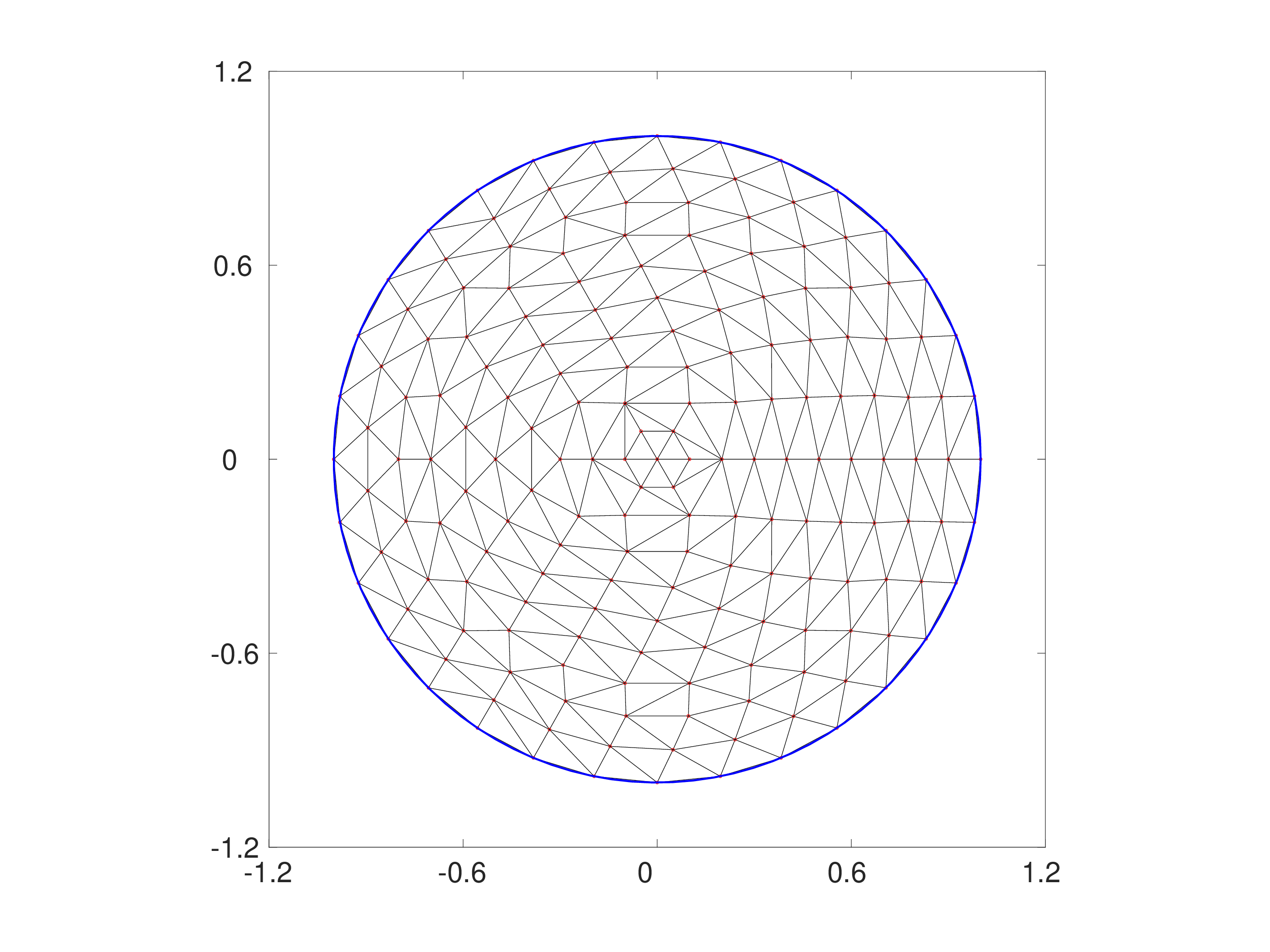}
     \caption{Triangulation: $180$ nodes, \\ $\mathrel{\phantom{==========}}$\; $316$ triangles. }
    \end{subfigure}
    \begin{subfigure}[t]{0.46\linewidth}
    \centering
    \includegraphics[scale=0.14]{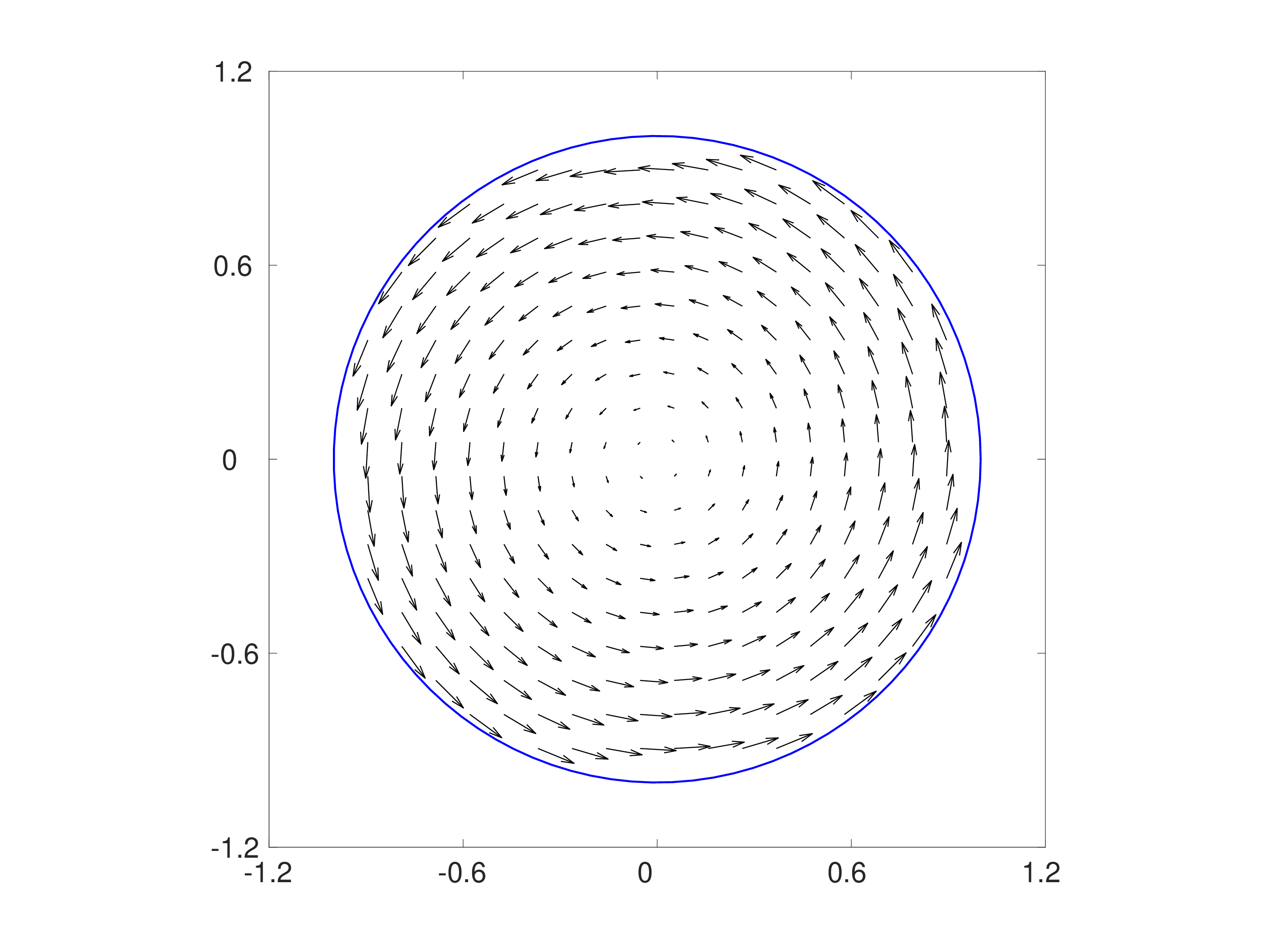}
   \caption{Body force: $f$}
    \end{subfigure}
   \caption{Initial configuration for the FEM implementation.}
   \label{fig: stokes-initial}
\end{figure}

\begin{figure}[htb!]
    \centering
     \begin{subfigure}[t]{0.46\linewidth}
    \centering
    \includegraphics[scale=0.14]{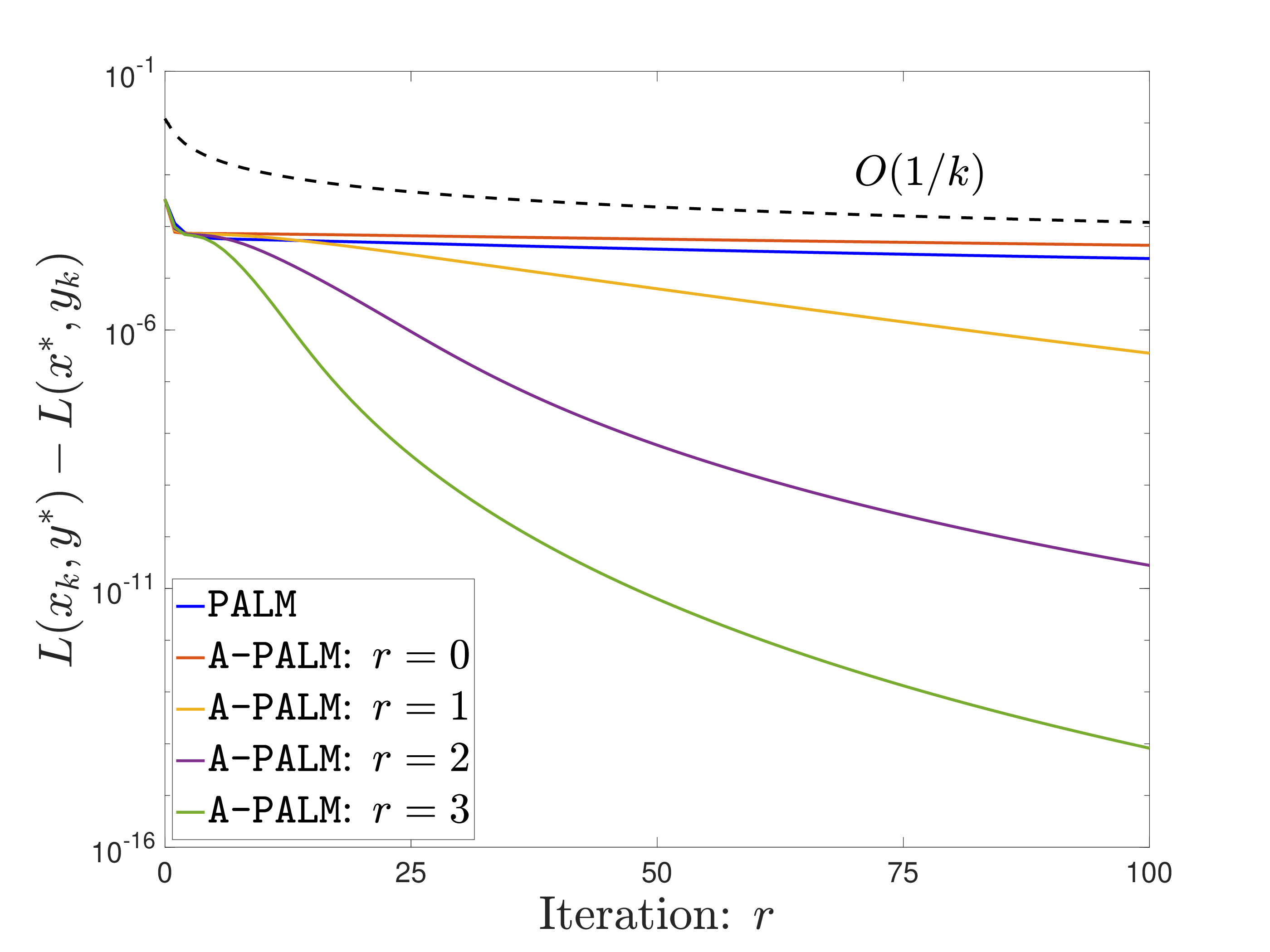}
     \caption{Primal-dual objective gap}
    \end{subfigure}
    \begin{subfigure}[t]{0.46\linewidth}
    \centering
    \includegraphics[scale=0.14]{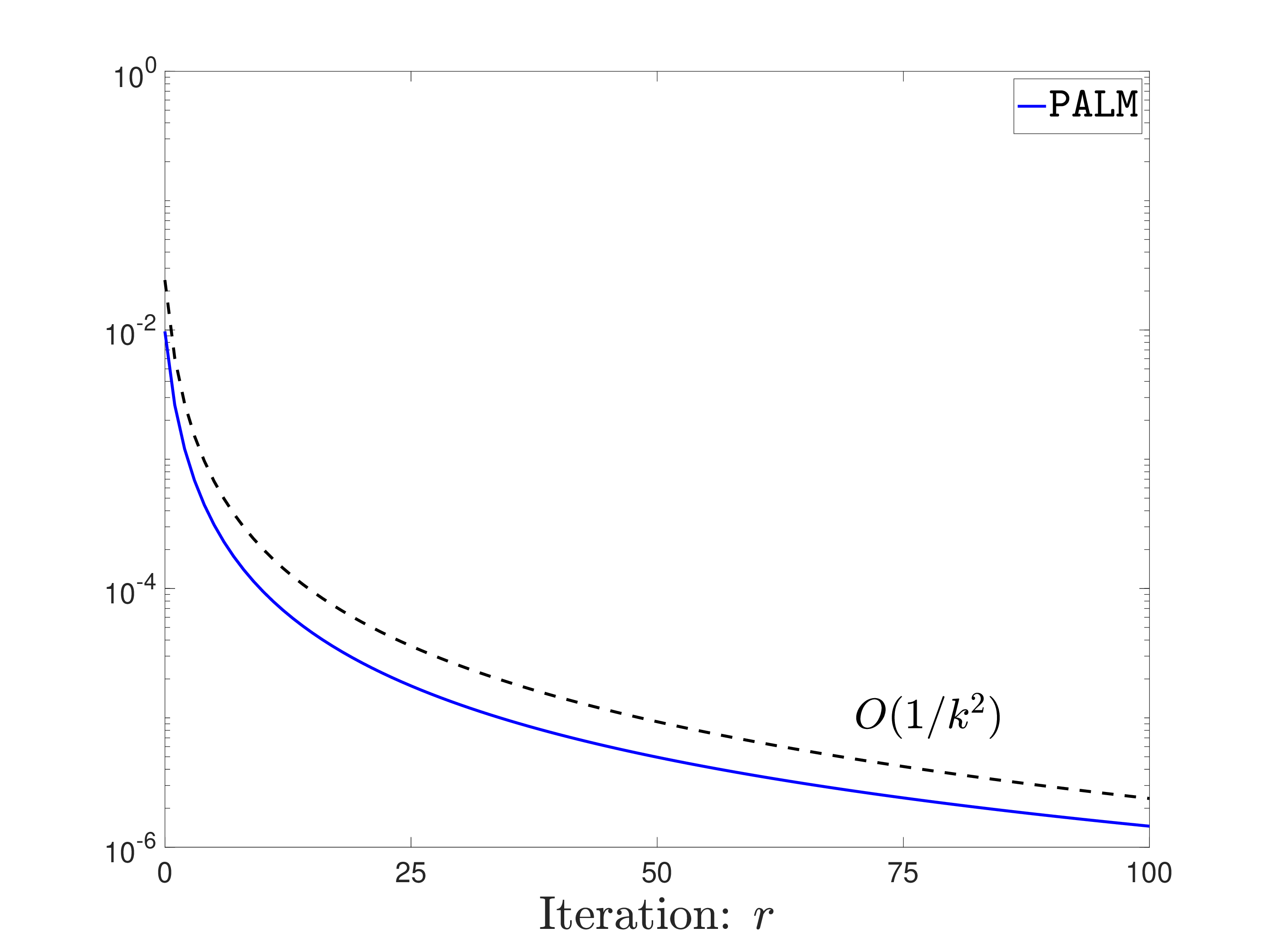}
   \caption{Gradient norm square}
    \end{subfigure}
   \caption{Convergence behavior of~\texttt{PALM} and~\texttt{A-PALM} on the steady Stokes system. Panel (a) shows the primal-dual objective gap for both algorithms, with~\texttt{A-PALM} set for parameters $r  \in \{0, 1, 2, 3\}$. Panel (b) shows the squared gradient norm for \texttt{PALM}}. 
   \label{fig: stokes}
\end{figure}

\begin{figure}[htb!]
    \centering
    \includegraphics[scale=0.16]{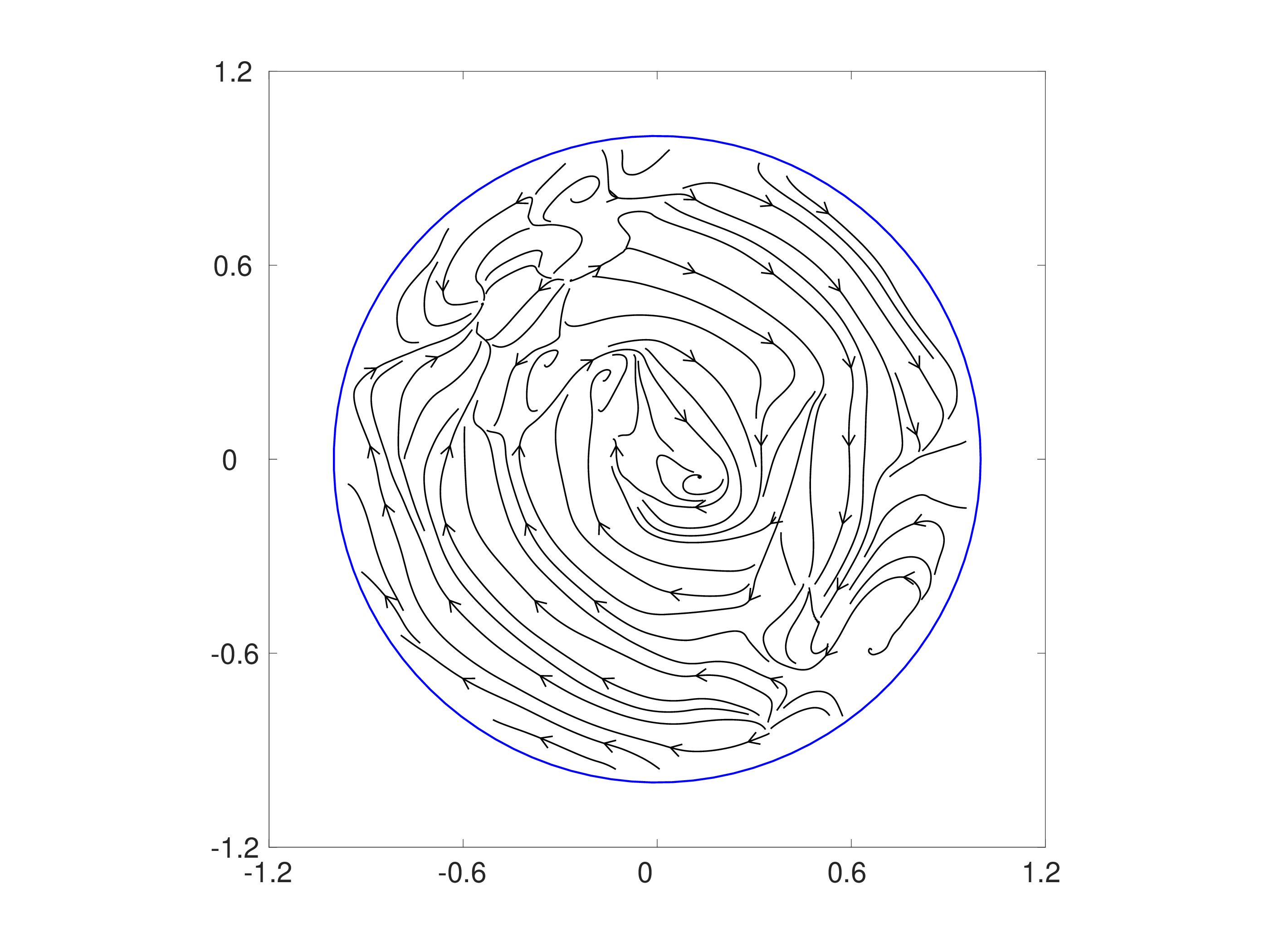}
   \caption{Numerical solution represented via the stream function $\Psi$}
   \label{fig: stokes-soln}
\end{figure}

To solve the system numerically, we employ the finite element method (FEM) for spatial discretization using quadratic elements to ensure sufficient accuracy for the velocity field.~\Cref{fig: stokes-initial} illustrates the initial triangular mesh used for the domain $\Omega$, as well as the body force field $\mathbf{f} = (-x_2, x_1)$ applied in the FEM implementation.
The numerical performance of~\texttt{PALM} and~\texttt{A-PALM} with parameters $r  \in \{0, 1, 2, 3\}$ is demonstrated in~\Cref{fig: stokes}.  Consistent with the benchmark test, \texttt{A-PALM} demonstrates significant acceleration as $r$ increases: for $r=3$, the primal-dual objective gap reaches $10^{-13}$ within 100 iterations. 
Furthermore, the primal-dual objective gap decay, coupled with the observed convergence of the squared gradient norm, confirms that our implementation aligns with the theoretical convergence rates established in~\Cref{thm: im-gad-dis-rate} and~\Cref{thm: im-n-dis-rate}. Finally, we present the numerical solution for the steady-state Stokes system~\eqref{eqn: stokes} obtained via~\texttt{A-PALM} with $r=3$ at iteration $k=100$. Given a reference point $(x_{1,0}, x_{2,0})$, the stream function $\Psi$ is defined as:
\[
\Psi = \int^{(x_1, x_2)}_{(x_{1,0}, x_{2,0})} -u_2 d\xi_1 + u_1 d\xi_2
\]
The resulting flow field, visualized via the stream function, is shown in~\Cref{fig: stokes-soln}. 

%
%
%
%
\section{Conclusion and further directions}
\label{sec: conclusion}

In this study, we demonstrate that the outer iteration of the proximal~\texttt{ALM}s is equivalent to the implicit~\texttt{GDA} schemes, which can also be naturally extended to accommodate the variable step sizes. Regarding the Lyapunov analysis, we highlight a fundamental distinction between convex and minimax problems: the potential function must transition from utilizing the objective gap to employing the variational inequality. Motivated by these two perspectives, we propose the variable step-size implicit~\texttt{GDA} scheme~\eqref{eqn: im-gda-discrete}, the family of Nesterov-type accelerated implicit GDA schemes~\eqref{eqn: im-n-gda}, and the Nesterov-type explicit~\texttt{GDA} scheme~\eqref{eqn: explicit-n-gda}. We develop the Lyapunov analysis for these methods to establish their convergence rates via the continuous-time ODE frameworks. For the implicit (proximal) schemes, the inner subproblems can be obtained efficiently using classical numerical solvers, such as successive over-relaxation or, more generally, the conjugate gradient method. Thus, accelerating the outer iteration is of significant interest in scientific computing, with widespread applications such as solving elliptic systems that model fluid flow (e.g., groundwater) through porous media~\citep{chen1999augmented}. Furthermore, the proposed gradient-based optimization methods require neither specialized preconditioners nor the large memory resources often associated with classical Krylov subspace methods~\citep{benzi2005numerical}.

Building on these two perspectives, we identify the application of Lyapunov analysis to gradient-based optimization methods for nonsmooth optimization, particularly proximal algorithms such as the alternating direction method of multipliers (\texttt{ADMM})~\citep{li2024understanding} and  the primal-dual hybrid gradient (\texttt{PDHG}) method~\citep{li2024understanding1, zeng2025lyapunov}, emerges as a compelling avenue for future research. Since nonsmooth problems generally preclude the use of classical Krylov subspace methods~\citep{benzi2005numerical}, accelerating the outer iterations thus becomes especially important.  Given the implicit nature of both \texttt{ADMM} and \texttt{PDHG}, constructing new Lyapunov functions through the lens of  the implicit \texttt{GDA} framework provides a very natural approach with the potential to recover classical ergodic convergence rates. Although the design of corresponding accelerated schemes remains highly challenging, this perspective opens a promising avenue for advancing future research in nonsmooth optimization. 


\section*{Acknowledgements}
This work was partially supported by SIMIS (startup fund and cross-disciplinary research projects) and by the NSFC (Grant No. 12241105). 

\bibliographystyle{abbrvnat}
\bibliography{sigproc}

\end{document}